\documentclass[a4,10pt,twoside]{article}
\usepackage{amsmath, amssymb, amsthm, amsfonts, amsxtra, latexsym, amscd,
pb-diagram, graphics}
\theoremstyle{plain}
\newtheorem{dl}{Theorem}[section]

\newtheorem{md}[dl]{Proposition}

\newtheorem{dn}[dl]{Definition}

\begin{document}
\renewcommand{\proofname}{$\mathbf{Proof}$}
\renewcommand{\refname}{$\mathbf{References}$}
\pagestyle{headings}

\title{\large{ \textbf{ STRUCTURE OF ANN-CATEGORIES AND\\ MAC LANE-SHUKLA COHOMOLOGY}}}
\author{\textbf{Nguyen Tien Quang}}

\date{}
\maketitle
\pagestyle{myheadings} 
\markboth{Structure of Ann-categories and MacLane - Shukla cohomology }{Nguyen Tien Quang }
\maketitle
\setcounter{tocdepth}{1}

{\bf Abstract.} In this paper\footnote{This paper has been published  in East-West J. of Mathematics Vol. 5, No 1, 2003} we study the structure of a class of categories having two operations which satisfy axioms      analoguos to that of rings. Such categories are called "Ann-categories". We obtain the classification theorems for      regular Ann-categories and Ann-functors by using Mac Lane-Shukla cohomology of rings. These      results give new interpretations of the cohomology groups $H^3(R, M)$ and $H^2(R, M)$ of the rings      $R$.

\section{\small{ Introduction and Preliminaries}}\label{S:P*}
$\!$\indent{Monoidal categories and symmetric monoidal categories were studied  first by S. Mac Lane [8], J. B\'{e}nabou [1] and G. M. Kelly [3]. They are, respectively, categories $\mathcal{A}$ together with a bifunctor $\otimes$: $\mathcal{A}\times\mathcal{A}\rightarrow\mathcal{A}$ and a system of natural equivalences of associativity-unitivity, or a system of natural equivalences of associativity-unitivity-commutativity. A. Solian [14], H. X. Sinh [2] and K. H. Ulbrich [15] , investigated $\otimes$-categories from   the point of view of algebraic structure. They examined the monoidal categories whose all objects are invertible.}\\
\indent{The problem of coherence always plays a fundamental role in the study of any class of $\otimes$-categories. From initial conditions, we have to prove that the morphisms generated by a given ones depend only on its source butt. The consideration of structures arose later in the papers of H. X. Sinh [2] and B. Mitchell [9]. Here we obtained deep results on the classification by the cohomology of groups.}\\
\indent{By the other direction, M. Laplaza [4] considered the coherence of natural equivalences of distributivity in a category having two operations $\oplus$ and $\otimes$. In the papers of Laplaza, the distribution of \emph{monomorphisms} together with the natural isomorphisms of the two symmetrical monoidal structures must satisfy 24 commutative diagrams, that form natural relations between them.\\
In this paper, we  consider a class of Pic-categories (see H. X. Sinh [2]) in which the second operation and natural equivalences of distributivity are defined so that the analogous axioms of rings are verified. Such categories are called Ann-categories. Coherence for Ann-categories was shown in [11].}\\
\indent{Throughout we define invariants of Ann-category basing on construction of reduced Ann-categories and pre-sticked of the type $(R, M)$. From this we obtain classification theorems for the regular Ann-categories and Ann-functors by using cohomology groups $H^3(R, M)$, $H^2(R, M)$ of the ring $R$. These theorems give a relation between the notion of Ann-category with the theory of cohomology of rings and the problem of extention of rings.}\\
\indent{For convinience, the tensor product of two objects $A$ and $B$ is denoted by $AB$ instead of $A\otimes B$, but for the morphisms we still write $f\otimes g$ to avoid confusion with the composition of morphisms.}\\
\indent{The notions and results on monoidal categories are supposed to be familier to the readers (see [3, 5, 8] for example).}\\ 

\indent{Recall that \emph{a Pic-category} is a symmetric monoidal category $\mathcal{A}$ (or a $\otimes$ACU-category $\mathcal{A}$) in which every object is invertible and every morphism is an isomorphism (see [2]).}\\
\begin{dn}
An Ann-category is a category $\mathcal{A}$ together with
\begin{enumerate}
\item[(i)] Two bifunctors $\;\oplus$, $\otimes$: $\mathcal{A}\times\mathcal{A}\rightarrow\mathcal{A}$.
\item[(ii)] A fixed object $0\in\mathcal{A}$ with natural isomorphisms $a^+, c, g, d$ such that $\big{(}\mathcal{A}, \oplus, a^+, c, (0, g, d)\big{)}$ is a Pic-category.
\item[(iii)] A fixed object $1\in\mathcal{A}$ with natural isomorphisms $a, l, r$ such that\\ $\big{(}\mathcal{A}, \otimes, a, (1, l,r)\big{)}$ is a monoidal category (i. e. a$\quad\otimes$AU-category).
\item[(iv)] Two natural isomorphisms $\mathfrak{L}$, $\mathfrak{R}$
$$\mathfrak{L}_{A,X,Y}:A(X\oplus Y)\rightarrow AX\oplus AY$$
$$\mathfrak{R}_{X,Y,A}:(X\oplus Y)A\rightarrow XA\oplus YA$$
satisfying the following conditions 
     \begin{enumerate}
         \item[(Ann-1)]For every object $A\in\mathcal{A}$, the pair of $\oplus$-functors $(L^A, \breve{L}^A)$, $(R^A, \breve{R}^A)$ defined by 
                   \begin{equation*}
                         \begin{cases}
                          L^A:X\rightarrow AX\\ 
                          \breve{L}^A_{X,Y}=\mathfrak{L}_{A,X,Y}
                         \end{cases}
                         \begin{cases}
                          R^A:X\rightarrow AX\\ 
                          \breve{R}^A_{X,Y}=\mathfrak{R}_{X,Y,A}
                         \end{cases}
                    \end{equation*} 
are $\oplus$AC-functors.
          \item[(Ann-2)]For any $A,B,X,Y\in\mathcal{A}$ the following diagrams are commutative\\\

\setlength\unitlength{0.5cm}
$$\begin{picture}(20,5)(3,0)
\put(0,0){$(AB)(X\oplus Y)$}
\put(5.5,0){\vector(4,0){12.8}}
\put(12,0.5){$\mathfrak{L}$}
\put(19,0){$(AB)X\oplus (AB)Y$}
\put(3,3){\vector(0,-1){2}}
\put(2.5,2){$a$}
\put(0,3.5){$A(B(X\oplus Y))$}
\put(5.5,3.5){\vector(4,0){3.2}}
\put(6,4){$id\otimes\mathfrak{L}$}
\put(9.5,3.5){$A(BX\oplus BY)$}
\put(15,3.5){\vector(4,0){3.2}}
\put(16.5,4){$\mathfrak{L}$}
\put(19,3.5){$A(BX)\oplus A(BY)$}
\put(22.2,3){\vector(0,-1){2}}
\put(22.7,2){$a\oplus a$}
\end{picture}$$

$$\begin{picture}(20,5)(3,0)
\put(0,0){$((X\oplus Y)A)B$}
\put(5.5,0){\vector(4,0){3.3}}
\put(6,0.5){$\mathfrak{R}\otimes id$}
\put(9.5,0){$(XA\oplus YA)B$}
\put(15,0){\vector(4,0){3.3}}
\put(16,0.5){$\mathfrak{R}$}
\put(19,0){$(XA)B\oplus (YA)B$}
\put(3,3){\vector(0,-1){2}}
\put(2,2){$a$}
\put(0,3.5){$(X\oplus Y)(AB)$}
\put(5.5,3.5){\vector(1,0){13}}
\put(11.5,4){$\mathfrak{R}$}
\put(19,3.5){$X(AB)\oplus Y(AB)$}
\put(22.2,3){\vector(0,-1){2}}
\put(22.7,2){$a\oplus a$}$$
\end{picture}$$   
 
$$\begin{picture}(20,5)(3,0)
\put(0,0){$(A(X\oplus Y))B$}
\put(5.5,0){\vector(4,0){3.3}}
\put(6,0.5){$\mathfrak{L}\otimes id $}
\put(9.5,0){$(AX\oplus AY)B$}
\put(15,0){\vector(4,0){3.3}}
\put(16,0.5){$\mathfrak{R}$}
\put(19,0){$(AX)B\oplus (AY)B$}
\put(3,3){\vector(0,-1){2}}
\put(2,2){$a$}
\put(0,3.5){$A((X\oplus Y)B)$}
\put(5.5,3.5){\vector(4,0){3.3}}
\put(6,4){$id\otimes\mathfrak{R}$}
\put(9.5,3.5){$A(XB\oplus YB)$}
\put(15,3.5){\vector(4,0){3.3}}
\put(16,4){$\mathfrak{L}$}
\put(19,3.5){$A(XB)\oplus A(YB)$}
\put(22.2,3){\vector(0,-1){2}}
\put(22.7,2){$a\oplus a$}
\end{picture}$$  

$$\begin{picture}(20,5)(3,0)
\put(0,0){$(AX\oplus BX)\oplus (AY\oplus BY)$}
\put(10.2,0){\vector(4,0){6}}
\put(12.9,0.5){v}
\put(16.8,0){$(AX\oplus AY)\oplus (BX\oplus BY)$}
\put(3,3){\vector(0,-1){2}}
\put(3.5,2){$\mathfrak{R}\oplus \mathfrak{R}$}

\put(0,3.5){$(A\oplus B)X\oplus (A\oplus B)Y$}
\put(9.5,3.5){\vector(-1,0){1}}
\put(8.7,4){$\mathfrak{L}$}
\put(10,3.5){$(A\oplus B)(X\oplus Y)$}
\put(16.5,3.5){\vector(4,0){1}}
\put(16.7,4){$\mathfrak{R}$}
\put(18,3.5){$A(X\oplus Y)\oplus B(X\oplus Y)$}
\put(22.2,3){\vector(0,-1){2}}
\put(22.7,2){$\mathfrak{L}\oplus \mathfrak{L}$}
\end{picture}$$
\noindent where $v=v_{A, B, C, D}$: $(A\oplus B)\oplus(C\oplus D)\rightarrow (A\oplus C)\oplus(B\oplus D)$ is the unique functorial morphism constructed from the morphisms $a^+, c$ and $id$ in the Pic-category $(\mathcal{A}, \oplus)$.
             \item[(Ann-3)]The following diagrams are commutative\\

$$\begin{CD}
1(X\oplus Y)@>\mathfrak{L}>>1X\oplus 1Y\\
@V\ell VV @VV\ell\oplus\ell V\\
X\oplus Y @= X\oplus Y
\end{CD}
\qquad\qquad
\begin{CD}
(X\oplus Y)1@>\mathfrak{R}>>X1\oplus Y1\\
@VrVV@VVr\oplus rV\\
X\oplus Y@=X\oplus Y
\end{CD}$$

     \end{enumerate}
\end{enumerate}
\end{dn} 
\begin{dn}
Let $\mathcal{A}$ and $\mathcal{A}'$ be Ann-categories . An Ann-functor from $\mathcal{A}$ to $\mathcal{A}^{'}$ is a functor $F$: $\mathcal{A}\rightarrow \mathcal{A}^{'}$ together with natural isomorphisms $\breve{F}$, $\widetilde{F}$ such that: $(F, \breve{F})$ is a $\oplus$AC-functor, $(F,\widetilde{F})$ is a $\otimes$A-functor and $\breve{F}$, $\widetilde{F}$ are compatible with natural equivalences of distributivity in the sense that the  following two diagrams are commutative\\
\setlength\unitlength{0.5cm}
\begin{picture}(20,6)(-1.2,0)
\put(0,0){$F(AB\oplus AC)$}
\put(5.5,0){\vector(4,0){3.3}}
\put(6,0.5){$\breve{F} $}
\put(9.1,0){$F(AB)\oplus F(AC)$}
\put(15.5,0){\vector(4,0){3.3}}
\put(16,0.5){$\widetilde{F}\oplus\widetilde{F}$}
\put(19,0){$FAFB\oplus FAFC$}
\put(3,3){\vector(0,-1){2}}
\put(0.5,2){$F(\mathfrak{L})$}
\put(0,3.5){$F(A(B\oplus C))$}
\put(5.5,3.5){\vector(4,0){3.3}}
\put(6,4){$\widetilde{F}$}
\put(9.5,3.5){$FAF(B\oplus C)$}
\put(15.2,3.5){\vector(4,0){3.5}}
\put(16,4){$id\otimes\breve{L}$}
\put(19,3.5){$FA(FB\oplus FC)$}
\put(23,3){\vector(0,-1){2}}
\put(23.5,2){$\mathfrak{L}^{'}$}
\end{picture}\\   
\begin{picture}(20,6)(-1.2,0)
\put(0,0){$F(AC\oplus BC)$}
\put(5.5,0){\vector(4,0){3.3}}
\put(6,0.5){$\breve{F} $}
\put(9.1,0){$F(AC)\oplus F(BC)$}
\put(15.5,0){\vector(4,0){3.3}}
\put(16,0.5){$\widetilde{F}\oplus\widetilde{F}$}
\put(19,0){$FAFC\oplus FBFC$}
\put(3,3){\vector(0,-1){2}}
\put(0.5,2){$F(\mathfrak{R})$}
\put(0,3.5){$F((A\oplus B)C)$}
\put(5.5,3.5){\vector(4,0){3.3}}
\put(6,4){$\widetilde{F}$}
\put(9.5,3.5){$F(A\oplus B)FC$}
\put(15.2,3.5){\vector(4,0){3.5}}
\put(16,4){$\breve{L}\otimes id$}
\put(19,3.5){$(FA\oplus FB)FC$}
\put(23,3){\vector(0,-1){2}}
\put(23.5,2){$\mathfrak{R}^{'}$}
\end{picture}\\   

\noindent If $F$ is an equivalence, then $(F,\breve{F},\widetilde{F})$ is called an Ann-equivalence. 
\end{dn} 
\begin{md}
Let $\mathcal{A}$ be an Ann-category and $A\in\mathcal{A}$. Then there exist unique isomorphisms $\widehat{L}^A$: $A\otimes 0\rightarrow 0$, $\widehat{R}^A$: $0\otimes A\rightarrow 0$ so that $(L^A,\breve{L}^A,\widehat{L}^A)$ and $(R^A,\breve{R}^A,\widehat{R}^A)$ are symmetrical monoidal functors ($\oplus$ACU-functor)
\end{md}
\begin{proof}Since $(\mathcal{A},\oplus)$ is a Pic-category, each $\oplus$AC-functor is also a $\oplus$ACU-functor.
\end{proof}
\begin{md}In any Ann-category $\mathcal{A}$, the isomorphisms $\widehat{L}^A$, $\widehat{R}^A$ have the following properties:
\begin{itemize}
\item[(i)]The family $\widehat{L}^-=\widehat{L}$ (resp. the family $\widehat{R}^-=\widehat{R}$) is a $\oplus$-morphism from the functor $(R^0, \breve{R}^0)$ (resp. $(L^0, \breve{L}^0)$) to the functor $(\theta:A\mapsto 0, \breve{\theta}=g_0^{-1})$ i. e. the following diagrams are commutative:
$$
\begin{CD}
A0@>f\otimes id>>B0\\
@V\widehat{L}^AVV@VV\widehat{L}^BV\\
0@=0
\end{CD}
\qquad\qquad
\begin{CD}
(X\oplus Y)0@>\breve{R}^0>>X0\oplus Y0\\
@VV\widehat{L}^{X\oplus Y}V @VV\widehat{L}^{X}\oplus\widehat{L}^{Y}V\\
0@<g_0<<0\oplus 0
\end{CD}
$$
(resp. $\widehat{R}^B(id\otimes f)=\widehat{R}^A$ and $\widehat{R}^{X\oplus Y}=g_0(\widehat{R}^X\oplus\widehat{R}^Y)\breve{L}^0)$.
\item[(ii)]For any $A, B\in\mathcal{A}$, the following diagrams are commutative:\\
\setlength\unitlength{0.5cm}
\begin{picture}(8,6)(-2,0)
\put(0,0){X0}
\put(1.5,0){\vector(1,0){2}}
\put(2,0.5){$\widehat{L}^X$}
\put(4,0){0}
\put(7,0){\vector(-1,0){2}}
\put(5.5,0.5){$\widehat{R}^Y$}
\put(7.5,0){0Y}
\put(1,3){\vector(0,-1){2}}
\put(-2.5,1.5){$id\otimes \widehat{R}^X$}
\put(0,3.5){X(0Y)}
\put(2.8,3.5){\vector(1,0){4}}
\put(5,4){a}
\put(7.5,3.5){(X0)Y}
\put(8,3){\vector(0,-1){2}}
\put(8.5,1.5){$\widehat{L}^X\otimes id$}
\put(13,0){(XY)0}
\put(15.7,0){\vector(1,0){2}}
\put(16.2,0.5){$\widehat{L}^{XY}$}
\put(18.2,0){0}
\put(13.9,3){\vector(0,-1){2}}
\put(13,1.5){a}
\put(13,3.5){X(Y0)}
\put(15.7,3.5){\vector(1,0){2}}
\put(15.6,4){$id\otimes\widehat{L}^Y$}
\put(18.2,3.5){X0}
\put(18.5,3){\vector(0,-1){2}}
\put(19,1.5){$\widehat{L}^X$}
\end{picture}\\\\
and $\widehat{R}^{XY}=\widehat{R}^Y(\widehat{R}^X\otimes id)a_{0,X,Y}$.
\item[(iii)] $L^1=l_0$, $R^1=r_0$.
\end{itemize}
\end{md}
\section{\small{The  first two invariants of an Ann-category}}\label{S:P*}
$\!$\indent{Let $\mathcal{A}$ be an Ann-category. Then the set $\Pi_0(\mathcal{A})$ of the isomorphic classes of objects of $\mathcal{A}$ is a ring with the operations induced by the ones $\oplus$, $\otimes$ in $\mathcal{A}$, and $\Pi_1(\mathcal{A})=Aut(0)$ is an abelian group with operation denoted by $+$.}\\
\indent{The following two Theorems on the structure of the Ann-categories can be found in [12].}\\
\begin{dl}
$\Pi_1(\mathcal{A})$ is an $\Pi_0(\mathcal{A})$-bimodule where the left and right operations of the ring $\Pi_0(\mathcal{A})$ on the abeian group $\Pi_1(\mathcal{A})$ are defined respectively by $$su=\lambda_X(u),\quad us=\rho_X(u),\quad X\in s\in\Pi_0(\mathcal{A}),\;u\in\Pi_1(\mathcal{A})$$ in which $\lambda_X,\;\rho_X$ are the two maps $Aut(0)\rightarrow Aut(0)$ given by the following commutative diagrams:
$$
\begin{CD}
X0@>\widehat{L}^X>>0\\
@Vid\otimes uVV @VV\lambda_X(u)V\\
X0@>\widehat{L}^X>>0
\end{CD}
\qquad\qquad\qquad
\begin{CD}
0X@>\widehat{R}^X>>0\\
@Vu\otimes idVV @VV\rho_{X}(u)V\\
0X@>\widehat{R}^X>>0
\end{CD}
$$
the following theorem shows the invariableness of $\Pi_0(\mathcal{A})$-bimodule $\Pi_1(\mathcal{A})$.
\end{dl} 
\begin{dl}
Given two Ann-categories $\mathcal{A}$, $\mathcal{A}^{'}$. Then any Ann-functor $(F,\breve{F},\widetilde{F})$: $\mathcal{A}\rightarrow\mathcal{A}^{'}$ yields a ring homomorphism
\begin{align*}
F_0:\quad \Pi_0(\mathcal{A})&\rightarrow \Pi_0(\mathcal{A}^{'})\\
clX&\mapsto clFX\\
\intertext{and a group homomorphism}
F_1:\quad \Pi_1(\mathcal{A})&\rightarrow \Pi_1(\mathcal{A}^{'})\\
u&\mapsto \gamma_{F0}^{-1}(Fu)
\end{align*}
having the properties $$F_1(su)=F_1(s)F_0(u)\qquad F_1(us)=F_0(u)F_1(s)$$ where $\gamma_A$: $Aut(0)\rightarrow Aut(0)$ is defined by $\gamma_A(u)=g_A(u\otimes id_A)g_A^{-1}$. Moreover, $F$ is an Ann-equivalence if and only if $F_0$, $F_1$ are isomorphisms. 
\end{dl}
\indent Hence $\Pi_0(\mathcal{A})$ and $\Pi_1(\mathcal{A})$ are the first two invariants of an Ann-category.
\section{\small{ Reduced Ann-categories}}
$\!$\indent{In preparing to define the third invariant of Ann-categories, we construct reduced Ann-categories. Let $\mathcal{A}$ be an Ann-category. The reduced category $\mathcal{S}$ in constructed from $\Pi_0(\mathcal{A})$ and $\Pi_1(\mathcal{A})$ as follows: its objects are the elements of $\Pi_0(\mathcal{A})$, its morphisms are the automorphisms of the form $(s,u)$ with \\$s\in\Pi_0(\mathcal{A}),\;u\in\Pi_1(\mathcal{A})$ i. e.}$$Aut(s)=\{s\}\times\Pi_1(\mathcal{A})$$ \indent{The composition law of morphisms is reduced by addition in $\Pi_1(\mathcal{A})$. We shall use the transmission of structures (see [10]) to change $\mathcal{S}$ into an Ann-category which is equivalent to $\mathcal{A}$. Choose for every $s\in\Pi_0(\mathcal{A})$ a representant $X_s\in\mathcal{A}$ such that $X_0=0,\, X_1=1$ and then, for every pair $s,\;t\in\Pi_0(\mathcal{A})$, two families of isomorphisms}$$\varphi_{s,t}:\;X_s\oplus X_t\rightarrow X_{s+t},\qquad\psi_{s,t}:\;X_sX_t\rightarrow X_{st}$$ such that $$\varphi_{0,t}=g_{X_t},\qquad\varphi_{s,0}=d_{X_s}$$
$$\psi_{1,t}=1_{X_t},\qquad\psi_{s,1}=r_{X_s},\qquad\psi_{0,t}=\widehat{R}^{X_t},\qquad\psi_{0,s}=\widehat{L}^{X_s}$$  
\indent{Defining the functor $H:\;\mathcal{S}\rightarrow\mathcal{A}$ by $H(s)=X_s$, $H(s,u)=\gamma_{X_s}(u)$ and putting $\breve{H}=\varphi^{-1},\;\widetilde{H}=\psi^{-1}$ we can use the theorem of transmission of structures (see [10]) to obtain $\mathcal{S}$ to be an Ann-category with the two operations in the explicit forms:}
\begin{align}
&s\oplus t=s+t\qquad (\text{sum in ring}\;\Pi_0(\mathcal{A})) \\
&(s,u)\oplus (t,v)=(s+t,u+v)\\
&s\otimes t=st\qquad (\text{product in ring}\;\Pi_0(\mathcal{A}))\\
&(s,u)\otimes (t,v)=(st,sv+ut)
\end{align} 
and with the natural equivalences induced by that of $\mathcal{A}$. $\mathcal{S}$ is called the reduced Ann-category of $\mathcal{A}$. We now have:
\begin{dl}
In the reduced Ann-category $\mathcal{S}$ of $\mathcal{A}$, the natural equivalences of unitivity of the two operations $\oplus$, $\otimes$ are identities, and the natural equivalences $\xi,\;\eta,\;\alpha,\;\lambda,\;\rho$ induced from $a^+,\;c,\;a,\;\mathfrak{L},\;\mathfrak{R}$ by $(H,\breve{H},\widetilde{H})$ are functions having the values in $\Pi_1(\mathcal{A})$ and satisfying the following relations
\begin{enumerate}
\item[1.]$\xi(y,z,t)-\xi(x+y,z,t)+\xi(x,y+z,t)-\xi(x,y,z+t)+\xi(x,y,z)=0$
\item[2.]$\xi(0,y,z)=\xi(x,0,t)=\xi(x,y,0)=0$
\item[3.]$\xi(x,y,z)-\xi(x,z,y)+\xi(z,x,y)-\eta(x,z)+\eta(x+y,z)-\eta(y,z)=0$
\item[4.]$\eta(x,y)+\eta(y,x)=0$
\item[5.]$x\eta(y,z)-\eta(xy,xz)=\lambda(x,y,z)-\lambda(x,z,y)$
\item[6.]$\eta(x,y)z-\eta(xz,yz)=\rho(x,y,z)-\rho(y,x,z)$
\item[7.]
  $x\xi(y,z,t)-\xi(xy,xz,xt)=\\
  \lambda(x,z,t)-\lambda(x,y+z,t)+\lambda(x,y,z+t)-\lambda(x,y,z)$
\item[8.]
  $\xi(x,y,z)t-\xi(xt,yt,zt)=\\
  \rho(y,z,t)-\rho(x+y,z,t)+\rho(x,y+z,t)-\rho(x,y,t)$
\item[9.]
  $\rho(x,y,z+t)-\rho(x,y,z)-\rho(x,y,t)+\lambda(x,z,t)\\
  +\lambda(y,z,t)-\lambda(x+y,z,t)=-\xi(xz+xt,yz,yt)\\
 +\xi(xz,xt,yz)-\eta(xt,yz)+\xi(xz+yz,xt,yt)-\xi(xz,yz,xt)$
\item[10.]
  $\alpha(x,y,z+t)-\alpha(x,y,z)-\alpha(x,y,t)=\\
  x\lambda(y,z,t)+\lambda(x,yz,yt)-\lambda(xy,z,t)$
\item[11.]
  $\alpha(x,y+z,t)-\alpha(x,y,t)-\alpha(x,z,t)=\\
  x\rho(y,z,t)-\rho(xy,xz,t)+\lambda(x,yt,zt)-\lambda(x,y,z)t$
\item[12.]
  $\alpha(x+y,z,t)-\alpha(x,z,t)-\alpha(y,z,t)=\\
-\rho(x,y,z)t-\rho(xz,yz,t)+\rho(x,y,zt)$
\item[13.]
  $x\alpha(y,z,t)-\alpha(xy,z,t)+\alpha(x,yz,t)\\
  -\alpha(x,y,zt)+\alpha(x,y,z)t=0$
\item[14.]$\alpha(1,y,z)=\alpha(x,1,z)=\alpha(x,y,1)=0$
\item[15.]$\alpha(0,y,z)=\alpha(x,0,t)=\alpha(x,y,0)=0$
\item[16.]$\lambda(1,y,z)=\lambda(0,y,,z)=\lambda(x,0,z)=\lambda(x,y,0)=0$
\item[17.]$\rho(x,y,1)=\rho(0,y,z)=\rho(x,0,z)=\rho(x,y,0)=0$
\end{enumerate} 
for $x,y,z,t\in\Pi_0(A)$.
\end{dl}
For the two choices of different representants $(X_s,\varphi,\psi)$, we can prove the followings:
\begin{md}
If $\mathcal{S}$ with $(X_s,\varphi,\psi)$ and $\mathcal{S}^{'}$ with $(X_s^{'},\varphi^{'},\psi^{'})$ are two reduced Ann-categories of $\mathcal{A}$, then there exists an Ann-equivalence $(F,\breve{F},\widetilde{F})$: $\mathcal{S}\rightarrow\mathcal{S}^{'}$, with $F=id$.
\end{md}
If we substitute $\Pi_0(\mathcal{A})$ by a ring $R$ and $\Pi_1(\mathcal{A})$ by an $R$-bimodule $M$, we can construct an Ann-category $\mathcal{I}$ with the operations $\oplus$, $\otimes$ defined by the relations (3.1)-(3.4) and the natural equivalences 
$$a^+=\xi,c=\eta,a=\alpha,\mathfrak{L}=\lambda,\mathfrak{R}=\rho$$ 
satisfying the relations in the theorem 3.1. This Ann-category $\mathcal{I}$ is called an \emph{Ann-category of type $(R,M)$}.\\
\indent{If the function $\eta$ satisfies the \emph{regular condition} $\eta(x,x)=0$, the family $(\xi,\eta,\alpha,\lambda,\rho)$ is a 3-cocycle of the ring $R$ with coefficients in the $R$-bimodule $M$ in the Mac Lane-Shukla sense (see theorem 4.3). In particular, when $\lambda=0,\rho=0,\xi=0$ we have $\eta=0$ and hence $\alpha$ becomes a normal 3-cocycle of the $\mathbb{Z}$-algebra $R$ in the Hochshild sense (see [10]).}\\
\indent{Any ring $R$ with the unit $1\not=0$ may be considered as an Ann-category of the type $(R,0)$. Hence we have proved the following theorem:}
\begin{dl}
Any Ann-category is an Ann-equivalence to an Ann-category of the type $(R,M)$.
\end{dl}
\section{\small{Cohomology classification of the regular Ann-categories}}
$\!$\indent{According to theorem 3.3 we have only to consider the classification of the Ann-categories having the first two common invariants.}
\begin{dn}
Let $R$ be a ring with unit, $M$ be an $R- bimodule$ considered as a ring with the null multiplication. An Ann-category $\mathcal{A}$ is called having pre-stick of the type $(R,M)$ if there exists a pair of ring isomorphisms $(\epsilon_0, \epsilon_1)$ 
$$\epsilon_0: R \longrightarrow \Pi_0(\mathcal{A}),\quad \epsilon_1: M \longrightarrow \Pi_1(\mathcal{A})$$
satisfying the conditions:
$$\epsilon_1(su) = \epsilon_0(s)\epsilon_1(u), \quad \epsilon_1(us) = \epsilon_1(u)\epsilon_0(s), \quad s \in R, u \in M.$$
\indent A morphism between two Ann-categories $\mathcal{A}, \mathcal{A'}$ having the same pre-stick of the type $(R, M)$ is an Ann-functor $(F,  \Breve{F},\widetilde{F}): \mathcal{A} \longrightarrow \mathcal{A'}$ such that the following diagrams  are commutative
$$
\begin {CD}
\Pi_0(\mathcal{A}) @> F_0 >> \Pi_0(\mathcal{A'})\\
@V \epsilon_0 VV @VV \epsilon_0 V\\
R @= R
\end{CD}
\qquad\qquad
\begin{CD}
\Pi_1(\mathcal{A}) @> F_1 >> \Pi_1(\mathcal{A'})\\
@V {\epsilon'}_1 VV @VV {\epsilon'}_1 V\\
M @= M
\end{CD}
$$
in which $F_0, F_1$ are two ring morphisms induced from $(F,  \Breve{F},\widetilde{F})$. It follows directly from  the definition that $F$ is an equivalence.\\
\indent The two Ann-categories $\mathcal{A}, \mathcal{A'}$ are called \emph{congruences} if there exists a morphism $(F,  \Breve{F},\widetilde{F})  between  them $.
 \end{dn}
\begin{dn}
An Ann-category $\mathcal{A}$ having a natural equivalence $c$ of commutativity so that $c_{X,X} = id$ is called a \emph{regular} Ann-category.
 \end{dn}
\indent For the regular Ann-categories we can define its third invariant, that is an element of Mac Lane - Shukla cohomology group $H^3(R,M)$ of the ring $R$.\\
\indent Recall that the cohomology of an algebra $\Lambda$ with coefficients in an \\ $\Lambda-$bimodule coincides with the Mac Lane cohomology of the ring $R = \Lambda$, considered as a $\mathbb{Z}-$algebra. We have
$$H^*(R, M) = H^*(\sum_{n \geq 0}Hom_\mathbb{Z}(U^n,M))$$ 
where $U$ is a graded differential algebra and a free resolution over $\mathbb{Z}$ of $R$. The differential $\delta$ over graded module $\sum Hom_\mathbb{Z}(U^n,M)$ is defined by the relation $\delta f= g + h$, where
$$g(u_1, ...,u_n)  = -\sum_{i=1}^{n}(-1)^{e_{i-1}}f(u_1, ..., du_i, ..., u_n),$$
\begin{multline*}
h(u_1, ...,u_n) = u_1f(u_2, ..., u_{n+1}) +\\
\sum_{i=1}^{n}(-1)^{e_i}f(u_1, \ldots, u_iu_{i+1},\ldots,u_{n+1}) + (-1)^{e_{n+1}}f(u_1,\ldots,u_n)u_{n+1},
\end{multline*}
$e_0 = 0, e_i = i + deg u_1 + \cdots + deg u_i \quad $ (see [13]).
\begin{dl}
A 3-cochain $f= <\zeta, \eta, \alpha,\lambda,\rho>$ of the ring $R$ with coefficients in the $R-$bimodule $M$ is a 3-cocycle if and only if $(\zeta, \eta, \alpha,-\lambda,\rho)$  is a family of natural equivalences of a regular Ann-category of the type $(R,M)$. 
\end{dl}
\noindent \begin{proof}
The essence of the proof is to compute the group $\mathbb{Z}^3(R,M)$ by choosing a convenient resolution of the ring $R$ (as a $\mathbb{Z}$-algebra), different from the two resolutions of Shukla and Mac Lane. For the additional structure of $R$, we consider the complex of abelian groups:
$$0 \longrightarrow B_4  \stackrel{d_4}{\longrightarrow}  B_3  \stackrel{d_3}{\longrightarrow}  B_2   \stackrel{d_2}{\longrightarrow}  B_1  \stackrel{d_1}{\longrightarrow}  B_0  \stackrel{\nu}{\longrightarrow}   R  \longrightarrow 0$$
in which
$$B_0 = \mathbb{Z}(\Dot{R}), \quad B_1 = \mathbb{Z}(\Dot{R}\times\Dot{R}), \quad B_2 = \mathbb{Z}(\Dot{R}\times\Dot{R}\times\Dot{R}) \oplus \mathbb{Z}(\Dot{R}\times\Dot{R})$$
$$B_3 = \mathbb{Z}(\Dot{R}\times \Dot{R}\times \Dot{R}\times \Dot{R}) \oplus \mathbb{Z}(\Dot{R}\times\Dot{R}\times\Dot{R})\oplus \mathbb{Z}(\Dot{R}\times\Dot{R}) \oplus \mathbb{Z}(\Dot{R})$$
$$B_4 =  Ker d_3, \quad  \Dot{R} = R \backslash   \{0\}$$
($\mathbb{Z}({\Dot{R}}^i),\, i = 1, 2, 3, 4$ are the free abelian groups generated by the set ${\Dot{R}}^i$).\\
\indent The morphisms  are given by:
\begin{eqnarray*}
\nu[x] & = & x, \quad x \in \Dot{R}\\
d_1[x,y] & = & [y] - [x+y] + [x] \\
d_2[x,y,z]  &=&[x,z] - [x+y,z] + [x,y+z] -[x,y]\\
d_2[x,y] &=& [x,y] - [y,x]\\
d_3[x,y,z,t] &=& [y,z,t] - [x+y,z,t] + [x, y+z,t] - [x,y,z+t] + [x,y,z]\\
d_3[x,y,z] &=& [x,y,z] - [x,z,y] + [z,x,y] + [x+y,z] - [x,z] - [y,z]\\
d_3[x,y] &=& [x,y] + [y,x] \\
d_3[x] &=& [x,x] 
\end{eqnarray*}
$d_4 = i$ is the natural embedding.\\
\indent We now define a distributive multiplication in $B = \sum B_i$ such that $B$ becomes a graded differential algebra over $\mathbb{Z}$. A 3-cochain $f$ is an element of a direct sum 
\begin{multline*}
Hom_\mathbb{Z}(B_2,M) \oplus Hom_\mathbb{Z} (B_1\otimes B_0, M) \oplus Hom_\mathbb{Z}(B_0\otimes B_1,M) \\
\oplus Hom_\mathbb{Z}(B_0 \otimes B_0\otimes B_0,M)$$
\end{multline*}
\indent This implies that $f$ is defined by a family of mappings 
\begin{eqnarray*}
\zeta(x,y,z) &=&f([x,y,z])\\
\eta(x,y) &=&f([x,y])\\
\lambda(x,y,z) &=& f([x] \otimes [y,z])\\
\rho (x,y,z)&=&f([x,y] \otimes [z])\\
\alpha(x, y,z) &=& f([x] \otimes [y] \otimes [z])
\end{eqnarray*} 
\indent From the formula of differentiation of the above resolution we complete the proof.
\end{proof}
\begin{dl}[Classification theorem]
There exists a bijection between the set of the congruence classes of pre-sticked regular Ann-categories of the type $(R,M)$ and the cohomology group $H^3(R,M)$ of the ring $R$, with coefficients in the $R$-bimudule $M$.
\end{dl}
\begin{proof}
Consider the resolution that is shown in the proof of the theorem 4.3. If $f = <\zeta,\eta,\alpha,\lambda,\rho>$ is 3-coboundary, $f = \delta g$, with $g$ is a pair of mappings 
\begin{eqnarray*}
\mu &:&B_1 \longrightarrow M\\
\nu &:&B_0 \otimes B_0 \longrightarrow M
\end{eqnarray*}
\noindent we have the following relations 
\begin{eqnarray*}
-\zeta(x,y,z) & = & \mu(y,z) - \mu(x+y,z) + \mu(x,y+z) - \mu(x,y)\\
-\eta(x,y) & = & \mu(x,y) - \mu(y,x) = ant \mu(x,y)\\
\alpha(x,y,z) & = & x\nu(y,z)- \nu(xy,z) + \nu(x,yz)-\nu(x,y)z\\
-\lambda(x,y,z) &=& \nu(x,y+z) - \nu(x,y) - \nu(x,z) + x\mu(y,z) - \mu(xy,xz)\\
\rho(x,y,z)& =& \nu(x+y,z) - \nu(x,z) - \nu(y,z) -\mu(x,y)z + \mu(xz,yz)
\end{eqnarray*}
These relations imply what we have to prove.
\end{proof}
$\!$\indent This theorem leads to the investigation of application of the Ann-category concept into the theory of ring extensions. The classtification theorem  in the general case is still an open problem. 
\section{\small{Ann-functors and low dimension cohomology groups of rings}}
$\!$\indent In this section given problem is that of finding whether there is Ann-functor between two Ann-categories and, if so, how many. Since each Ann-category is Ann-equivalent to one Ann-category of the type $(R,M)$ so the solution of problem for a class of Ann-categories of the type $(R,M)$ is enough.\\
\indent If $f = <\zeta,\eta,\alpha,\lambda,\rho>$ is a 3-cocycle in $\mathbb{Z}^3(R,M)$ the structure \\ $(\zeta,\eta,\alpha,-\lambda,\rho)$ of Ann-category $(R,M)$ is denoted by $\Hat{f}$. Moreover, if
$$F = (F,  \Breve{F},\widetilde{F}) : (R,M,\Hat{f}) \longrightarrow (R',M',\Hat{f'})$$
is an Ann-functor, this functor is a pair of ring homomorphisms $(F_0,F_1)$ compatible with actions of bimodule. So sometimes $F$ is denoted by $(F_0,F_1)$. $R'$-bimodule $M$ may be changed into $R$-bimodule by the homomorphism $F_0$,
$$m'r = mF(r), \quad rm' = F(r)m', \quad r \in R, m' \in M'.$$
Because $f \in \mathbb{Z}^3(R,M) $ and $ f' \in \mathbb{Z}^3(R',M'), \quad F$ induces canonically 3-cocycles
$$f_{*}, {f'}^{*} \in \mathbb{Z}^3(R,M').$$ 
For axample 
\begin{eqnarray*}
\zeta_{*}(x,y,z) = F(\zeta(x,y,z))\\
{\zeta '}^*(x,y,z) = \zeta '(Fx,Fy,Fz).
\end{eqnarray*}
Isomorphisms $\Breve{F},\widetilde{F}$ are mappings $R\times R \longrightarrow M'$
\[ \begin{array}{ c c c c c c  c}
\mu (x,y)&=& \Breve{F}_{x,y}&:&F(x+y) &\longrightarrow& Fx + Fy\\
\nu(x,y)&=&\widetilde{F}_{x,y}&:&F(xy) &\longrightarrow & (Fx)(Fy)
\end{array} \]
These mappings, according to definition, satisfy diagrams in definition 1.2. On the other hand, $<\mu,\nu>$ is a 2-cochain of ring cohomology. From a calculation of $H^3(R,M)$ we have
\begin{equation}
f_* - {f'}^* =\delta<\mu,\nu>
\end{equation}
\begin{dl}
Let  $\mathcal{I}= (R,M, \Hat{f}),\;\mathcal{I'} = (R', M',\hat{f'}) $ be two regular Ann-categories and 
$$F= (F_0,F_1): \mathcal{I} \longrightarrow\mathcal{I'} $$
be a functor that satisfies the condition (5.1). Then $F$ is an Ann-functor if and only if $H_*(f) - H^*(f') = 0$ in $H^3(R,M')$. In this case, we can say that Ann-functor $(F,  \Breve{F},\widetilde{F})$ is induced by the functor  $F$.
\end{dl}
\noindent \begin{proof}
If $(F,  \Breve{F},\widetilde{F})$ is an Ann-functor with $ \Breve{F} = \mu, \widetilde{F}=\nu$, the condition (5.1) gives equation
$$H_*(f) - H^*(f') = 0$$
\indent Conversely, the equation $H_*(f) - H^*(f') = 0$ automatically implies $f_* - {f'}^* = \delta g$, there $g=<\mu,\nu>$ is a 2-cochain. Let $\Breve{F} = \mu, \widetilde{F} = \nu$, we have an Ann-functor $(F,  \Breve{F},\widetilde{F})$.
\end{proof}
\begin{dn}
An Ann-functor $F:(R,M,f) \longrightarrow (R',M',f')$ is called regular if $F$ satisfies condition $f_* = {f'}^*$.
\end{dn}
\indent In case there exists a regular Ann-functor $F$, we have the following theorem
\begin{dl}
(i) There exists a bijection between the set of the congruence classes of regular Ann-functors induced by a pair $(F_0,F_1)$ and the cohomology group $H^2(R,M')$ of the ring $R$ with coefficients in the $R$-bimodule $M'$.\\
\indent (ii) If $F: (R,M,f) \longrightarrow (R',M',f')$ is an Ann-functor, there exists a bijection $$Aut(F) \longrightarrow \mathbb{Z}^1(R,M')$$ between the group of automorphisms of Ann-functor $F$ and the group $\mathbb{Z}^1(R,M')$.
\end{dl}
\noindent \begin{proof}
$(i)$ Let  $(F,  \Breve{F},\widetilde{F})$ be a regular Ann-functor
$$(F,  \Breve{F},\widetilde{F}) : (R,M,f) \longrightarrow (R',M',f)$$
Then $$f_* - {f'}^* = \delta<\mu,\nu> = 0$$ where $ \Breve{F} = \mu, \widetilde{F}=\nu$. It means $<\mu,\nu>$ is 2-cocycle.\\
\indent Suppore that $(G,  \Breve{G},\widetilde{G})$ is another regular Ann-functor $$(G,  \Breve{G},\widetilde{G}) : (R,M,f) \longrightarrow (R',M',f)$$
and $\alpha : F \longrightarrow G$ is an Ann-morphism. Then, by to definition, the following 
diagrams are commutative
$$\begin {CD}
F(x+y) @> \Breve{F} >> Fx + Fy\\
@V \alpha_{x+y}VV @VV \alpha_x + \alpha_yV\\
G(x+y) @> \Breve{G} >> Gx + Gy
\end{CD}$$
$$\begin {CD}
F(xy) @> \widetilde{F} >> (Fx)( Fy)\\
@V \alpha_{xy}VV @VV \alpha_x \otimes \alpha_yV\\
G(xy) @> \widetilde{G} >> (Gx )(Gy)
\end{CD}$$
where $x,y\in R$. Also from the definition we have
$$\alpha_x \otimes \alpha_y = (Fx)\alpha_y + \alpha_x(Fy) = x \alpha_y + \alpha_x y$$
so
\begin{eqnarray*}
 \Breve{G}_{x,y} - \Breve{F}_{x,y} &= &\alpha_x - \alpha_{x+y} + \alpha_y\\
\widetilde{G}_{x,y} - \widetilde{F}_{x,y} &=& x\alpha_y - \alpha_{x+y} + \alpha_x y.
\end{eqnarray*}
Because $ g= <\Breve{F},\widetilde{F}>, \quad  g' =  <\Breve{G},\widetilde{G}>$ are 2-cocycles and $\alpha$ is 1-cochain and by a calculation of $H^2(R,M)$ we have
\begin{equation}
g' - g = \delta\alpha
\end{equation}
\indent Equation (6) proves the existance of a correspondence from a class of regular Ann-functors cls$(F,  \Breve{F},\widetilde{F})$ to a class of cohomologies $g+ B^2(R,M'), \\ g= <\Breve{F},\widetilde{F}>$. Moreover this correspondence is an injection. We now prove that it is a projection. In fact, let $g=<\mu,\nu>$ be any 2-cocycle. Then we can directly verify that $(F,\mu,\nu)$ is a regular Ann-functor $(R,M,f)$ to $(R',M',f)$ corresponding to 2-cocycle $g$, proving (i).\\
\indent $(ii)$ Let $$F=(F,\mu,\nu):(R,M,f) \longrightarrow (R',M',f)$$ be an Ann-functor and $\alpha \in Aut(F)$. Then the equation (6) becomes $\delta(\alpha) = 0$, i.e. $\alpha \in \mathbb{Z}^1(R,M')$, proving (ii).
\end{proof}
\section{\small{Ann-category and theory of the extensions of rings}}
$\!$\indent In this section, we establish a direct relation between theory of the extensions of rings and theory of Ann-categories. According to Mac Lane [7] we call a \emph{bimultiplication} of a ring $A$ a pair of mappings $a \longmapsto \sigma a, a \longmapsto a \sigma$ of $A$ into itself which satisfy the rules
\begin{alignat*}{5}
\sigma(a + b) &= & \sigma a + \sigma b \quad &,&\quad (a+b)\sigma &= & a \sigma + b \sigma\\
\sigma(ab)& = & (\sigma a)b \quad &,& \quad (ab)\sigma &=& a (b\sigma)\\
&&a(\sigma b) \quad &=& (a \sigma)b
\end{alignat*}
for all elements $a, b \in A$. The sum $\sigma + \nu$ and the product $\sigma\nu$ of two bimultiplications $\sigma$ and $\nu$ are defined by the equations
\begin{alignat*}{5}
(\sigma + \nu )a & = & \sigma a + \nu a \quad &, &\quad a(\sigma + \nu) & = & a \sigma + a \nu \\
(\sigma  \nu)a & = & \sigma(\nu a) \quad & , & \quad a(\sigma \nu) & = & (a \sigma)\nu
\end{alignat*}
for all $a$ in $A$.\\
\indent The set of all bimultiplications of $A$ is a ring denoted by $M_A.$ For each element c of $A$, a bimultiplication $\mu_c$ is defined by 
$$\mu_c a = ca,\quad a\mu_c = ac,\quad a \in A$$
\noindent We call $\mu_c$ an \emph{inner bimultiplication}. Clearly $\mu : A \longrightarrow M_A$ is a ring homomorphism and the image $\mu A$ of this homomorphism is a two-sided ideal in $M_A$. The quotient ring $P_A = M_A/\mu A$ is called the ring of \emph{outer bimultiplications} of $A$ and ring homomorphism $\theta : R \longrightarrow P_A$ is called \emph{regular} if $\theta (1) = 1$ and two any elements of $\theta (R)$ are \emph{permutable} ( the bimultiplications $\sigma$ and $\nu$ are called permutable if $\sigma (a \nu) = (\sigma a)\nu $ and $\nu (a \sigma) = (\nu a)\sigma $ for every $a$ in $A$). Then 
$$C_A = \{c \in A | ca = ac = 0 , \forall a \in A\} $$
is called \emph{bicenter} of $A$, and $C_A$ is a $R$-bimodule under the operations
$$xc = (\theta x)c ,\quad cx = c(\theta x) , \quad c \in C_A , x\in A.$$ 
\indent The "Extention problem" of rings requires finding the exact sequence of rings
$$0 \longrightarrow A\longrightarrow S \longrightarrow R \longrightarrow1  $$
induces homomorphism  $\theta : R \longrightarrow P_A$.\\
\indent Let $\sigma : R \longrightarrow M_A$ be a mapping such that $\sigma(x) \in \theta x, x \in R $ and $\sigma (0) = 0,$\\$ \sigma (1) = 1$. Then we define two mappings 
\begin{eqnarray*}
f : R \times R &\longrightarrow&A \\
g : R \times R &\longrightarrow&A 
\end{eqnarray*}
such that 
\begin{eqnarray*}
\mu f(x,y) &=&\sigma (x+y) - \sigma (x) - \sigma (x)\\
\mu g(x,y) &=&\sigma (xy) - \sigma (x) \sigma (x)
\end{eqnarray*}
for all $x,y \in R.$ The  ring structure of $M_A$ implies mappings $ \zeta ,\alpha ,\lambda ,\rho : {M_A}^3 \longrightarrow C_A$ and $\eta : {M_A}^2 \longrightarrow C_A$
\begin{eqnarray*}
\zeta(x,y,z) & = &f(x,y) - f(x+y,z) + f(x,y+z) - f(x,y)\\
\eta(x,y) & = &f(x,y) -f(y,x)\\
\alpha(x,y,z) & = &xg(y,z) -g(x,y,z) + g(x,y,z) - g(x,y)z \\
\lambda(x,y,z) & = &xf(y,z) - f(xy,xz) + g(x,y+z) - g(x,y) - g(x,z)\\
\rho(x,y,z) & = &f(x,y)z - f(xz,yz) + g(x+y,z) - g(x,z) -g(y,z)
\end{eqnarray*}
\indent We call the family ($\zeta ,\eta , \alpha , \lambda , \rho$) of the above mappings an obstruction of the regular homomorphism $\theta$. We can prove that if all these mappings are null, the homomorphism  $\theta : R \longrightarrow P_A$ can be realized by a ring extention. It is the ring 
$$S = \{ (a,r) \quad | \quad a \in A, r \in R \} $$
with operations
\begin{eqnarray*}
(a_1,r_1) +(a_2,r_2) &=& (a_1 +a_2 + f(r_1,r_2), r_1 +r_2 )\\
(a_1,r_1) (a_2,r_2) &=& (r_1a_2 +a_1r_2 + g(r_1r_2), r_1r_2 )
\end{eqnarray*}
\indent In the general case we have
\begin{md}
If ($\zeta ,\eta , \alpha , \lambda , \rho$) is an \emph{obstruction} of the regular homomorphism $\theta : R \longrightarrow P_A$, it is a family of natural equivalences of Ann-categories of the type $(R, C_A)$.
\end{md}
\noindent \begin{proof}
We can verify directly that $\zeta ,\eta , \alpha , \lambda , \rho$ satisfy the relations in the proposition 3.1 .
\end{proof}

Math. Dept., Hanoi University of Education\\
E-mail adresses: nguyenquang272002@gmail.com
\end{document}